\documentclass[12pt]{amsart}

\usepackage{amssymb}
\usepackage{latexsym}
\usepackage{amsmath}

\def\c{{\mathbf C}}
\def\l{{\mathcal L}}
\def\a{{\mathcal A}}
\def\x{{\mathcal X}}

\def\pp{{\mathbf P}} 
\def\uu{{\mathbf U}}
\def\ss{{\mathbf S}}
\def\oo{{\mathbf O}}
\def\nn{{\mathbf N}}
\def\ww{{\mathbf W}}

\begin{document}

\newtheorem{defi}{Definition}[section]
\newtheorem{prop}{Proposition}[section]
\newtheorem{theo}{Theorem}[section]
\newtheorem{lemm}{Lemma}[section]
\newtheorem{coro}{Corollary}[section]

\title
{Symmetries of a generic coaction}
\author{Teodor Banica}
\address{
Institut de Math\'ematiques de Jussieu, Case 191, Universit\'e Paris 6, 
4 Place Jussieu, 75005 Paris}
\email{banica@math.jussieu.fr}
\thanks{e-mail address: banica@math.jussieu.fr}

\begin{abstract}
If $B$ is $\c^*$-algebra of dimension $4\leq n<\infty$ then the finite
dimensional irreducible representations of the compact quantum
automorphism group of $B$, say
$G_{aut}(\widehat{B})$, have the same fusion rules as the ones of
$\ss\oo (3)$. As consequences, we get (1) a structure result for
$G_{aut}(\widehat{B})$ in the case where $B$ is a matrix algebra (2)
if $n\geq 5$ then the dual $\widehat{G}_{aut}(\widehat{B})$ is not
amenable (3) if $n\geq 4$ then the fixed point subfactor $P^{G_{aut}(\widehat{B})}\subset (B\otimes P)^{G_{aut}(\widehat{B})}$ has index $n$ and principal graph $A_\infty$.
\end{abstract}

\maketitle

\section*{Introduction}

Let $X_n$ be the space consisting of $n$ points. The category of
groups acting on $X_n$ has a universal object: the symmetric group
${\mathcal S}_n$. S. Wang has recently discovered that
when replacing ``groups'' with ``compact quantum groups'' the
resulting category has also a universal object, say $G_{aut}(X_n)$. If
$n=1,2,3$ then $G_{aut}(X_n)={\mathcal S}_n$. If $n\geq 4$ this compact quantum
group $G_{aut}(X_n)$ is not a classical group, nor a finite quantum group (see
  \cite{wang}, \cite{ergo}) and very less
  seems to be known about it. Some quantum subgroups of
  it, which are of interest in connection with spin models, were constructed in \cite{spin}. 

More generally, any finite dimensional $\c^*$-algebra $B$ has a
compact quantum group of automorphisms, say $G_{aut}(\widehat{B})$
(we have $G_{aut}(\widehat{\c}^n)=G_{aut}(X_n)$). See \cite{wang}. See also theorem 1.1 below
and the comments preceding it: actually when $B$ is noncommutative the definition of $G_{aut}(\widehat{B})$
requires as data a trace on $B$ -- i.e. one can define only compact
quantum groups of automorphisms of {\em measured} quantum finite spaces, cf. \cite{wang} -- and the distinguished trace we use
here is in general different from the distinguished trace used in \cite{wang}.

In this paper we prove that if $dim(B)\geq 4$ then the set of classes of finite dimensional irreducible representations of $G_{aut}(\widehat{B})$ can be labeled by the positive integers, $Irr(G_{aut}(\widehat{B}))=\{ p_n\mid n\in\nn\}$, such that the fusion rules
are
$$p_k\otimes p_s=p_{\mid k-s\mid}+p_{\mid k-s\mid
  +1}+\dots +p_{k+s-1}+p_{k+s}$$
In other words, we prove that we have an isomorphism
  of fusion semirings
$$R^+(G_{aut}(\widehat{B}))\simeq R^+(\ss\oo (3))$$
This kind of isomorphisms already appeared in quite various contexts,
and seem to be deeply related to notions of ``deformation''. See the
survey \cite{survey}.

The reasons for the existence of the above isomorphism are somehow
hidden by the technical details of the proof, and will be explained now.

``4'' comes from Jones' index (cf. proposition 2.2);  and also
from Wang's computations for $n=1,2,3$. The fact that the
distinguished trace on $B$ we use is the ``good'' one comes from
Markov inclusions (cf. proposition 2.1 (ii)); and also from our result
itself. For other traces proposition 2.1 (iv) shows that the
fundamental representation of the corresponding $G_{aut}(\widehat{B})$ contains,
besides the trivial representation, at least two components. This
situation reminds the one of $A_o(F)$'s for $F\bar{F}\notin \c\cdot Id$ (see \cite{iso} for what happens in this case).

The category $Rep(G_{aut}(\widehat{B}))$ of representations of $G_{aut}(\widehat{B})$ is in a
certain sense ``generated'' by two arrows: the multiplication $\mu
:B\otimes B\to B$ and the unit $\eta :\c\to B$ (see proposition 1.1). The
point is that the only ``relations'' satisfied by these ``generators''
are the ``universal'' ones coming from the axioms of the algebra structure
of $B$ (see lemma 2.1). This shows that $Rep(G_{aut}(\widehat{B}))$ ``does not depend so
much on $B$''. More precisely, with a good index and a good trace, one can show that its semiring of objects -- which is $R^+(G_{aut}(\widehat{B}))$ --
does not depend on $B$. On the other hand one can see from definitions (see also corollary 4.1) that for
$B=M_2(\c )$ we have
$$G_{aut}(\widehat{M_2(\c )})\simeq \pp\uu (2)\simeq\ss\oo (3)$$
with $\pp\uu (2)$ acting on $M_2(\c )$ in the obvious way. Thus
$R^+(G_{aut}(\widehat{B}))\simeq R^+(\ss\oo (3))$.

The paper is organized as follows. In section 1 we find
convenient ``presentations'' of the Hopf $\c^*$-algebra of continuous
functions on $G_{aut}(\widehat{B})$ and of its category of corepresentations. In
section 2 we use theory from \cite{j}, \cite{ghj} for finding the good
index and the good trace on $B$. In
sections 3 and 4 we use techniques from \cite{cras} for proving the
main result. By combining it with results from \cite{cras},
\cite{subf}, \cite{spin} we get the consequences (1), (2), (3) listed in the abstract.

\section{Coactions and corepresentations}

There is only one Hopf $\c^*$-algebra to be considered in this paper, namely
$A^{aut}(B)$, the object dual to the compact quantum group $G_{aut}(\widehat{B})$ in
the introduction. $A^{aut}(B)$ is by definition a certain $\c^*$-algebra
given with generators and relations. Its precise structure is that
of a finitely generated full Woronowicz-Kac algebra.

However, for understanding the definition of $A^{aut}(B)$ we have to
state one of its universality properties. The really relevant one
holds in the category of finitely generated full Woronowicz-Kac
algebras.

The definition of these algebras is as follows (see \cite{survey} for
explanations on terminology).

\begin{defi}[cf. definition 1.1 in \cite{w1}]
A finitely generated (or co-matricial) full Woronowicz-Kac algebra is
a pair $(A,u)$ consisting of a unital $\c^*$-algebra $A$ and a
unitary matrix $u\in M_n(A)$ satisfying the following conditions:

\noindent (i) $A$ is the enveloping $\c^*$-algebra of its $*$-subalgebra $\a$ generated by the entries of $u$.

\noindent (ii) there exists a $\c^*$-morphism $\Delta :A\rightarrow A\otimes A$ such that $(id\otimes\Delta )u=u_{12}u_{13}$.

\noindent (iii) there exists a $\c^*$-antimorphism $S:A\to A$ sending
$u_{ij}\leftrightarrow u_{ji}^*$.
\end{defi}

One can define a $\c^*$-morphism $\varepsilon :A\to\c$ by $\varepsilon
(u_{ij})=\delta_{i,j}$. The maps $\Delta ,\varepsilon ,S$ satisfy the
well-known requirements for a comultiplication, counit and
antipode. More precisely, their restrictions to $\a$ make $\a$ a Hopf
$\c$-algebra. See \cite{w1}. 

Let $V$ be a finite dimensional $\c$-linear space. A coaction of the
Hopf algebra $\a$ on $V$ is a linear map $\beta :V\to V\otimes
\a$ satisfying
$$(id\otimes\Delta )\beta =(\beta\otimes id )\beta ,\,\,\,\,
(id\otimes\varepsilon )\beta =id$$
A corepresentation of $\a$ on $V$ is an element $u\in\l (V)\otimes\a$
satisfying
$$(id\otimes\Delta )u=u_{12}u_{13},\,\,\,\, (id\otimes\varepsilon
)u=1$$
The coactions (resp. corepresentations) of the Hopf algebra $\a$ will
be called coactions (resp. corepresentations) of the finitely
generated full Woronowicz-Kac algebra $A$. These two notions are related as follows.

\begin{lemm}
Denote by $\beta\mapsto u_\beta$ the composition of canonical linear isomorphisms
$${\mathcal L} (V,V\otimes\a )\simeq V^*\otimes  V\otimes \a\simeq {\mathcal L}
(V)\otimes \a$$
and by $u\mapsto \beta_u$ its inverse. That is, if $\{v_i\}$
is a basis of $V$ and $\{ e_{ij} :v_j\mapsto v_i\}$ is the corresponding system of matrix
units in $\l (V)$, then the relation between $\beta$ and $u_\beta$ is
$$\beta (v_i)=\sum_j v_j\otimes u_{ji}\Longrightarrow
u_\beta =\sum_{ij}e_{ij}\otimes u_{ij}$$
If $\beta :V\to V\otimes \a$ is a linear map, then $\beta$ is a
coaction iff $u_\beta$ is a corepresentation.\hfill $\Box$
\end{lemm}

Let $B$ be a finite dimensional $\c^*$-algebra. We denote by $\mu
:B\otimes B\to B$ the multiplication and by $\eta :\c\to B$ the linear
map which sends $1\mapsto 1$. Let $tr:B\to\c$ be a faithful positive
normalised trace (in what follows we call such traces just ``traces''). We recall that the scalar product $<x,y>=tr(y^*x)$ makes $B$ into a Hilbert space. 

Let $A$ be a finitely generated full Woronowicz-Kac algebra. A coaction of $A$ on
$(B,tr)$ is a morphism of unital $*$-algebras $\beta :B\to B\otimes\a$ which is
a coaction of $A$ on the $\c$-linear space $B$, and which is such that
the trace satisfies the equivariance property
$$(id\otimes tr)\beta =tr(.)1$$

We will use the following notations. If $D$ is a unital $\c$-algebra and $V,W$
are two finite dimensional $\c$-linear spaces and $v\in\l (V)\otimes D$ and $w\in\l (W)\otimes D$ we define
$$v\otimes w=v_{13}w_{23}\in \l (V)\otimes\l (W)\otimes D$$
$$Hom(v,w)=\{ T\in \l (V,W)\mid (T\otimes id)v=w(T\otimes id)\}$$
If $D$ is a bialgebra and $v$ and $w$ are corepresentations, then
$\otimes$ and $Hom$ are the usual tensor product and space of
intertwiners. In general, it is possible to associate to any algebra
$D$ a certain monoidal category having these $Hom$ and $\otimes$, so
our notations are not as abusive as they seem.

\begin{lemm}
Let $u\in\l (B)\otimes\a$ be a corepresentation, and denote $\beta :=\beta_u$.

(i) $\beta$ is multiplicative $\Longleftrightarrow$ $\mu\in Hom(u^{\otimes
  2},u)$.

(ii)  $\beta$ is unital $\Longleftrightarrow$ $\eta\in Hom(1,u)$.

(iii) $tr$ is $\beta$-equivariant $\Longleftrightarrow$ $\eta\in
Hom(1,u^*)$.

If these conditions are satisfied, then:

(iv)  $\beta$ is involutive $\Longleftrightarrow$ $u$ is unitary.
\end{lemm}

\begin{proof}
Let $\{b_1,b_2,...,b_n\}$ be an orthonormal basis of
the Hilbert space $B$. We will use many times the formula 
$$x=\sum_i tr(b_i^*x)b_i=\sum_i tr(xb_i^*)b_i=\sum_i tr(b_ix)b_i^*=\sum_i tr(xb_i)b_i^*$$
for any $x\in B$, and especially its particular case $1=\sum
tr(b_i^*)b_i=\sum tr(b_i)b_i^*$. For any $i,j$ we have
$$(\mu\otimes 1)u^{\otimes 2}(b_i\otimes b_j\otimes 1)=
(\mu\otimes 1)(\sum_{kl}b_k\otimes b_l\otimes u_{ki}u_{lj})=
\sum_{kl}b_kb_l\otimes u_{ki}u_{lj}=\beta (b_i)\beta (b_j)$$
$$u(\mu\otimes 1)(b_i\otimes b_j\otimes 1)=
u(\sum_kb_k\otimes tr(b_ib_jb_k^*)1)=
\sum_{kl} b_l\otimes
tr(b_ib_jb_k^*)u_{lk}=\beta (b_ib_j)$$
and (i) follows. The assertion (ii) is clear from
$$u(1\otimes 1)=u(\sum_jb_j\otimes tr(b_j^*)1)=\sum_{ij}b_i\otimes
tr(b_j^*)u_{ij}=\sum_jtr(b_j^*)\beta (b_j)=\beta (1)$$
As for (iii), we have $1\otimes 1=\sum_j b_j\otimes (tr(b_j)1)^*$ and
$$u^*(1\otimes 1)=
(\sum_{ij}e_{ji}\otimes u_{ij}^*)(\sum_Ib_I\otimes tr(b_I^*)1)=
\sum_{ij}b_j\otimes
tr(b_i^*)u_{ij}^*=\sum_jb_j\otimes ((tr\otimes id)\beta (b_j))^*$$
Let us prove (iv). Assume that (i,ii,iii) are satisfied and that
$\beta$ is involutive. Then for any $i,k$ we have
$$\sum_j u_{ji}^*u_{jk}=(tr\otimes id)(\sum_{js} b_j^*b_s\otimes
u_{ji}^*u_{sk})=(tr\otimes id)\beta
(b_i^*b_k)=tr(b_i^*b_k)1=\delta_{i,k}$$
so $u^*u=1$. On the other hand, as $u$ is a corepresentation of the
Hopf algebra $\a$, we know that $u$ is invertible (its inverse is
$(id\otimes S)u$: this follows by considering $(id\otimes E)u$,
with $E=m(S\otimes id)\Delta =m(id\otimes S)\Delta =\varepsilon
(.)1$). Thus $u$ is a unitary. Conversely, assume that (i,ii,iii)
are satisfied and that $u$ is unitary. From $\mu\in Hom(u^{\otimes
  2},u)$ we get $\mu^*\in Hom(u,u^{\otimes
  2})$, and together with $\eta\in Hom(1,u)$ this gives $\mu^*\eta\in
Hom(1,u^{\otimes 2})$. As $u^{\otimes 2}=u_{13}u_{23}$ it follows that 
$$u_{23}(\mu^*\eta (1)\otimes 1)=u_{13}^*(\mu^*\eta (1)\otimes 1)$$
We have $<\mu^*\eta (1),b_p\otimes b_q>=<1,b_pb_q>=tr(b_p^*b_q^*)$ for
any $p,q$, so that
$$\mu^*\eta (1)=\sum_{pq}tr(b_p^*b_q^*)b_p\otimes b_q$$
Let us compute $u_{23}(\mu^*\eta (1)\otimes 1)$ and
$u_{13}^*(\mu^*\eta (1)\otimes 1)$ by using this formula:
$$u_{23}(\mu^*\eta (1)\otimes 1)=
(\sum_{is}id\otimes e_{is}\otimes
u_{is})(\sum_{jS}tr(b_j^*b_S^*)b_j\otimes b_S\otimes 1)=
\sum_{jis}b_j\otimes b_i\otimes tr(b_j^*b_s^*)u_{is}$$
$$u_{13}^*(\mu^*\eta (1)\otimes 1)=
(\sum_{jk}e_{jk}\otimes id\otimes
u_{kj}^*)(\sum_{iK}tr(b_K^*b_i^*)b_K\otimes b_i\otimes 1)=
\sum_{jik}b_j\otimes b_i\otimes
tr(b_k^*b_i^*)u_{kj}^*$$
Thus from $u_{23}(\mu^*\eta (1)\otimes 1)=u_{13}^*(\mu^*\eta (1)\otimes 1)$
we get that
$$\sum_{s}tr(b_j^*b_s^*)u_{is}=\sum_{k} tr(b_k^*b_i^*)u_{kj}^*$$
for any $i$. The fact that $\beta$ is involutive follows from this and from
$$\beta (b_j^*)=\sum_s\beta (b_s)tr(b_j^*b_s^*)=\sum_{is}b_i\otimes
tr(b_j^*b_s^*)u_{is}$$
$$\beta (b_j)^*=\sum_kb_k^*\otimes u_{kj}^*=\sum_{ik} b_i\otimes
tr(b_k^*b_i^*)u_{kj}^*$$
\end{proof}

We will use now lemmas 1.1 and 1.2 for associating to any pair $(B,tr)$ a
certain algebra $A^{aut}(B,tr)$ and a category ${{\mathcal C}_{B,tr}}$. The theorem
1.1 below claims no originality and may be found, in a slightly
different form, in \cite{wang}. Notice the following two differences
between it and theorem 5.1 in \cite{wang}.

(1) Only the case of a certain distinguished trace was explicitely
worked out in \cite{wang}, and the general case was left in there to the
reader. What happens is that Wang's trace is not the ``good'' one,
i.e. it is in general different from the one needed for having
``minimality'' of $R^+(A^{aut}(B,tr))$ (see the comments in the
introduction), which is called ``canonical
trace'' in section 2 below. By the way, this is the reason why we will
use later on the notation $A^{aut}(B)$ (for $A^{aut}(B,tr)$ with
$tr$= the canonical trace in the sense of definition 2.1 below) instead of Wang's notation $A_{aut}(B)$ (which
corresponds to $A^{aut}(B,tr)$ with
$tr$= the distinguished trace used in \cite{wang}): these two algebras may
not be isomorphic in general.

(2) The presentation of $A^{aut}(B,tr)$ given here -- to be used in
proposition 1.1 for finding a ``presentation'' of its category of corepresentations -- is different from the one in \cite{wang}. The point
is that this kind of presentation is the ``good'' one in a certain (quite
obvious) sense (see section 2 in \cite{survey}). Of course one
can prove, via manipulations of generators and relations, or just by
using uniqueness of universal objects, that our algebra
$A^{aut}(B,tr)$ is the same as Wang's $A_{aut}(B,tr)$. For reasons
of putting $aut$ as an exponent in our notation see the above comment (1).

\begin{theo}[cf. \cite{wang}]
Let $(B,tr)$ be a finite dimensional $\c^*$-algebra together with a
trace. Denote by $\mu :B\otimes B\to B$ the multiplication and by
$\eta :\c\to B$ the linear unital map. Choose an orthonormal basis of $B$, and use it for identifying $B\simeq \c^n$ as Hilbert
spaces, with $n=dim(B)$. Consider the following universal $\c^*$-algebra
$A^{aut}(B,tr)$:
$$\c^*<(u_{ij})_{i,j=1,...,n}\mid u=(u_{ij})\,\,\mbox{is
  unitary},\, \eta\in Hom(1,u),\, \mu\in Hom(u^{\otimes 2},u)>$$

(i) There exists a unique structure of finitely generated full
  Woronowicz-Kac algebra on $A^{aut}(B,tr)$ which
  makes $u$ a corepresentation. There exists a unique coaction $\beta$
  of $A^{aut}(B,tr)$ on $(B,tr)$ such that $u_\beta =u$ via the above identification $B\simeq \c^n$.

(ii) If $A$ is a finitely generated full Woronowicz-Kac algebra and
$\gamma$ is a coaction of $A$ on $(B,tr)$ then there exists a unique
morphism $f:A^{aut}(B,tr)\to A$ such that $(id\otimes f)\beta
=\gamma$. Moreover, $(A^{aut}(B,tr),\beta )$ is the unique pair
(finitely generated full Woronowicz-Kac algebra, coaction of it on
$B$) having this property.
\end{theo}

The definition of $A^{aut}(B,tr)$ should be understood as follows. Let
$F$ be the free $*$-algebra on $n^2$ variables $(u_{ij})_{i,j=1,...,n}$ and let
$u=(u_{ij})\in\l (\c^n)\otimes F$. By explicitating the notations for
$Hom$ and $\otimes$ with $D=F$ we see that both conditions $\eta\in
Hom(1,u)$ and $\mu\in Hom(u^{\otimes 2},u)$, as well as the
condition ``$u$ is unitary'', could be interpreted as being a
collection of relations between the $u_{ij}$'s and their adjoints. Let
$J\subset F$ be the two-sided $*$-ideal  
generated by all these relations. Then the matrix $u=(u_{ij})$ is unitary in $M_n({\bf C} )\otimes (F/J)$, so its coefficients $u_{ij}$ are of norm 
$\leq 1$ for every  ${\bf C}^*$-seminorm on $F/J$ and the enveloping
${\bf C}^*$-algebra of $F/J$ is well-defined. We call it
$A^{aut}(B,tr)$. The discussion on the (in)dependence of
$A^{aut}(B,tr)$ on the basis of $B$ is left to the reader.

We will use freely the terminology from \cite{w2} concerning concrete
monoidal $\ww^*$-categories. We recall that the word ``concrete''
comes from the fact that the monoidal $\ww^*$-category is given together with
an embedding into (= faithful monoidal $\ww^*$-functor to) the
category of finite dimensional Hilbert spaces.

\begin{prop}
The concrete monoidal $\ww^*$-category
$Corep(A^{aut}(B,tr))$ of finite dimensional unitary smooth
corepresentations of $A^{aut}(B,tr)$ is the completion in the sense of
\cite{w2} of the concrete monoidal $\ww^*$-category ${{\mathcal C}_{B,tr}}$ defined as follows:

- the monoid of objects of ${{\mathcal C}_{B,tr}}$ is $(\nn ,+)$.

- the Hilbert space associated to an object $m\in\nn$ is $B^{\otimes m}$.

- ${{\mathcal C}_{B,tr}}$ is the smallest concrete monoidal $\ww^*$-category containing
the arrows $\eta$, $\mu$.
\end{prop}

The definition of ${{\mathcal C}_{B,tr}}$ should be understood as follows: its arrows
are linear combinations of (composable) compositions of tensor
products of maps of the form $\eta$, $\mu$, $\eta^*$, $\mu^*$ and
$id_m:=$ identity of $B^{\otimes m}$. It is clear that ${{\mathcal C}_{B,tr}}$ is a
concrete monoidal $\ww^*$-category.

\begin{proof}
We will prove both results at the same time. Let us consider the
concrete monoidal $\ww^*$-category ${{\mathcal C}_{B,tr}}$ in proposition 1.1. It is clear that
the pair $(A^{aut}(B,tr),u)$ in theorem 1.1 is its universal
admissible pair in the sense of \cite{w2}.

We prove now that the object $1$ of ${{\mathcal C}_{B,tr}}$ is a complex
conjugation for itself in the sense of \cite{w2}. Let us define an
invertible antilinear map $j:B\to B$ by $j(b_q)=b_q^*$ for any
$q$. With the notations from page 39 in \cite{w2} we have
$$t_j(1)=\sum_p b_p\otimes j(b_p)=\sum_p b_p\otimes
b_p^*=\sum_{pq}tr(b_p^*b_q^*)b_p\otimes b_q$$
We have seen in proof of lemma 1.2 that $\mu^*\eta (1)$ is given by
the same formula, and it follows that $t_j=\mu^*\eta$. In particular
we get that $t_j\in Hom_{{\mathcal C}_{B,tr}} (0,2)$. By choosing as basis $X=\{
b_1,b_2,\ldots ,b_n\}$ of $B$ a complete system of matrix units we may
assume that $X=X^*$. It follows that $j=j^{-1}$, so with the notations
in \cite{w2} we get that $\bar{t}_j=t_{j^{-1}}^*=t_j^*$ is in
$Hom_{{\mathcal C}_{B,tr}} (2,0)$. Thus $1=\bar{1}$ in ${{\mathcal C}_{B,tr}}$ in the sense of
\cite{w2}. Notice that this shows also that $1=\bar{1}$ in ${{\mathcal C}_{B,tr}}$ in
the sense of \cite{lr}. Moreover, from
$$\mid\mid t_j(1)\mid\mid^2=\sum_{pq}<b_p\otimes b_p^*,b_q\otimes
b_q^*>=n$$
we get that the dimension $d_{{\mathcal C}_{B,tr}} (1)$ of $1$ in ${{\mathcal C}_{B,tr}}$ in the sense
of \cite{lr} is equal to $n$.

As $1=\bar{1}$ theorem 1.3 in \cite{w2} applies and shows that
$(A^{aut}(B,tr),u)$ is a finitely generated Woronowicz algebra
(i.e. a compact matrix pseudogroup, with the terminology in there)
whose concrete monoidal $\ww^*$-category of corepresentations is the
completion of ${{\mathcal C}_{B,tr}}$. Also $A^{aut}(B,tr)$ is full by definition, and
it is of Kac type because the quantum dimension of its fundamental
corepresentation $u$ is $d_{{\mathcal C}_{B,tr}} (1)=n$, hence is equal to its
classical dimension (see \cite{lr}, see also section 1 in \cite{subf}).

Summing up, we have proved both the first assertion in theorem 1.1 (i) and
proposition 1.1. The other assertions in theorem 1.1 follow from lemmas 1.1 and 1.2.
\end{proof}

\section{Good trace, good index}

Let $B$ be a finite dimensional $\c^*$-algebra. Let $n=dim(B)$. Denote
by $id_s$ the identity of $\l (B^{\otimes s})$ for any $s$. As in section 1, we
denote by $\mu :B\otimes B\to B$ the multiplication and by $\eta
:\c\to B$ the linear unital map. 

We recall that each faithful trace on $B$ makes it into a Hilbert
space, so in particular it gives rise to adjoints $\mu^*:B\to B\otimes B$ and $\eta^*:B\to\c$.

\begin{prop}
If $tr:B\to\c$ is a faithful normalised trace then the following are
equivalent:

(i) $tr$ is the restriction of the unique trace of $\l (B)$, via the
embedding $B\subset \l (B)$ given by the left regular representation.

(ii) $\c\subset B$ is a Markov inclusion in the sense of \cite{ghj}.

(iii) if $\phi :B\simeq \bigoplus_{\gamma =1}^sM_{m_\gamma}$ is a
decomposition of $B$ as a multimatrix algebra, then the weights
$\lambda_\gamma :=tr(\phi^{-1} (Id_{M_{m_{\gamma}}}))$ of $tr$ are given
by $\lambda_\gamma =n^{-1}m_\gamma^2$
for any $\gamma$.

(iv) $\mu\mu^*=n\cdot id$.
\end{prop}

\begin{proof}
The equivalence between (i), (ii) and (iii) is clear from
definitions. Let us prove that (iii) and (iv) are equivalent. We may identify $B$ with a multimatrix algebra $\bigoplus_{\gamma
  =1}^sM_{m_\gamma}$ as in (iii). Let $\{\lambda_\gamma\}$ be the weights
  of $tr$. Then $tr(e_{ij}^\gamma
  )=\delta_{i,j}m_\gamma^{-1}\lambda_\gamma$ for any $\gamma ,i,j$, so the set 
$$\{f_{ij}^\gamma :=m_\gamma^{1/2}\lambda_\gamma^{-1/2}e_{ij}^\gamma\mid\gamma =1,2,...,s,\, i,j=1,2,...,m_\gamma\}$$
is an orthonormal basis of $B$. Thus for any $\gamma ,i,j$ we have
$$\mu^*(e_{ij}^\gamma )=
\sum_{\delta\varepsilon klpq}f_{kl}^\delta\otimes 
f_{pq}^\varepsilon <\mu^*(e_{ij}^\gamma )
,f_{kl}^\delta\otimes  f_{pq}^\varepsilon >=
\sum_{\delta\varepsilon klpq}f_{kl}^\delta\otimes 
f_{pq}^\varepsilon tr(f_{kl}^\delta f_{pq}^\varepsilon e_{ji}^\gamma )=$$
$$\sum_{\delta\varepsilon klpq}m_\delta^{1/2}\lambda_\delta^{-1/2}e_{kl}^\delta\otimes 
m_\varepsilon^{1/2}\lambda_\varepsilon^{-1/2}e_{pq}^\varepsilon
tr(m_\delta^{1/2}\lambda_\delta^{-1/2}e_{kl}^\delta
m_\varepsilon^{1/2}\lambda_\varepsilon^{-1/2}e_{pq}^\varepsilon e_{ji}^\gamma )=m_\gamma\lambda_\gamma^{-1}\sum_le_{il}^\gamma\otimes
  e_{lj}^\gamma$$
Thus for any $\gamma ,i,j$ we have $\mu\mu^* (e_{ij}^\gamma )=m_\gamma\lambda_\gamma^{-1}\sum_l
e_{ij}^\gamma =m_\gamma^2\lambda_\gamma^{-1}e_{ij}^\gamma$, so (iii)
$\Longleftrightarrow$ (iv).
\end{proof}

\begin{defi}
The distinguished trace in proposition 2.1, say $\tau$, will be called
the canonical trace of $B$. The finitely generated full Woronowicz-Kac
algebra $A^{aut}(B,\tau )$ in theorem 1.1 will be denoted
$A^{aut}(B)$. The concrete monoidal $\ww^*$-category ${\mathcal C}_{B,\tau}$ in
proposition 1.1 wil be denoted ${{\mathcal C}_B}$.
\end{defi}

The arrows $\eta$, $\mu$, $\eta^*$, $\mu^*$ could be thought of as being
``generators'' of ${{\mathcal C}_B}$. In the next lemma we collect the
relevant ``relations'' satisfied by these arrows.

\begin{lemm}
(i) $\mu\mu^*=n\cdot id$ and $\eta^*\eta =id$.

(ii) $(\mu\otimes id)(id\otimes\mu^*)=(id\otimes\mu )(\mu^*\otimes
id)=\mu^*\mu$.

(iii) $\mu (\mu\otimes id)=\mu (id\otimes\mu )$.

(iv) $\mu (id\otimes\eta )=\mu (\eta\otimes id)=id$.
\end{lemm}

\begin{proof}
We have $\eta^*\eta (1)=<\eta^*\eta (1),1>=<\eta (1),\eta
(1)>=<1,1>=tr(1)=1$. Also the equality $\mu\mu^*=id$ was already
proved, and (iii) and (iv) are trivial, so it remains to prove
(ii). For, we may use an identification $B=\bigoplus_{\gamma
  =1}^sM_{m_\gamma}$ as in proof of proposition 2.1. By using the
formula of $\mu^*$ in there (with $\lambda_\gamma
=n^{-1}m_\gamma^2$ !) we get
$$(\mu\otimes id)(id\otimes\mu^*)(e_{pq}^\delta\otimes e_{ij}^\gamma )=
nm_\gamma^{-1}\sum_l(\mu\otimes id)(e_{pq}^\delta\otimes e_{il}^\gamma\otimes
  e_{lj}^\gamma )=nm_\gamma^{-1}\delta_{\gamma
    ,\delta}\delta_{q,i}\sum_le_{pl}^\gamma\otimes e_{lj}^\gamma$$
$$(id\otimes\mu)(\mu^*\otimes id)(e_{pq}^\delta\otimes e_{ij}^\gamma )=
nm_\delta^{-1}\sum_l(id\otimes\mu )(e_{pl}^\delta\otimes e_{lq}^\delta\otimes
  e_{ij}^\gamma )=nm_\delta^{-1}\delta_{\gamma
    ,\delta}\delta_{q,i}\sum_le_{pl}^\gamma\otimes e_{lj}^\gamma$$
$$\mu^*\mu (e_{pq}^\delta\otimes e_{ij}^\gamma )=\delta_{\gamma
    ,\delta}\delta_{q,i}\mu^*(e_{pj}^\gamma )=\delta_{\gamma
    ,\delta}\delta_{q,i} nm_\gamma^{-1}\sum_le_{pl}^\gamma\otimes
    e_{lj}^\gamma$$
for any $\gamma ,\delta ,i,j,p,q$, and this proves (ii).
\end{proof}

We recall that for $m\in {\mathbf N}$ and $\beta >0$ the $m$-th Temperley-Lieb algebra
$A_{\beta ,m}$ of index $\beta$ is defined with generators $e_1,e_2,...,e_{m-1}$ and Jones'
relations (see \cite{j}):

$e_i=e_i^*=e_i^2$ for any $i$.

$e_ie_j=e_je_i$ for any $i$ and $j$ with $\mid i-j\mid\geq 2$.

$\beta e_ie_je_i=e_i$ for $i$ and $j$ with $\mid i-j\mid =1$.

\begin{prop}
Assume that $n\geq 4$ and let $k\in\nn$. 

(i) $P=n^{-1}\mu^*\mu\in\l (B\otimes B)$ and
$Q=\eta\eta^*\in\l (B)$ are projections satisfying
$$n(Q\otimes id)P(Q\otimes id)=(Q\otimes id)$$
$$n(id\otimes Q)P(id\otimes Q)=(id\otimes Q)$$
$$nP(id\otimes Q)P=nP(Q\otimes id)P=P$$
$$(P\otimes id)(id\otimes P)=(id\otimes P)(P\otimes id)$$

(ii) There exists a faithful representation $\pi_k:A_{n,2k}\to \l (B^{\otimes k})$ which sends
$$e_{2s}\mapsto id_{s-1}\otimes P\otimes id_{k-s-1}\,\,\,\,
(s=1,2,...,k-1)$$
$$e_{2s+1}\mapsto id_s\otimes Q\otimes id_{k-s-1}\,\,\,\,
(s=0,1,...,k-1)$$

(iii) $dim(End_{{\mathcal C}_B} (k))\geq C_{2k}$ (the $2k$-th Catalan number).
\end{prop}

\begin{proof}
The assertions in (i) follow via easy computations from the formulas
in lemma 2.1. From (i) we get Jones' relations, hence a representation
$\pi_k$ as in (ii). This representation is nothing but the well-known
one coming by applying basic constructions to the Markov inclusion
$\c\subset B$ (cf. proposition 2.1 (ii)), so it is faithful. As $\mu$
and $\eta$ are arrows of ${{\mathcal C}_B}$, it follows that $P\in End_{{\mathcal C}_B} (2)$
and $Q\in End_{{\mathcal C}_B} (1)$. Thus (ii) gives a copy of $A_{n,2k}$ into
$End_{{\mathcal C}_B} (k)$ for any $k$. On the other hand, as the index $n\geq 4$
is generic, we have by \cite{j} that $dim(A_{n,2k})=C_{2k}$, and this
proves (iii).
\end{proof}

We will see in next section that $End_{{\mathcal C}_B} (k)\simeq A_{n,2k}$ for any $k$.

\section{Computation of $Hom_{{\mathcal C}_B} (0,k)$}

We fix an algebra $B$ of dimension $n\geq 4$ and we use the notations in section 2.

\begin{lemm}
Let ${{\mathcal C}_B}^+$ be the set of arrows of ${{\mathcal C}_B}$ consisting of linear
combinations of compositions of tensor products of maps
of the form $\eta$, $\mu^*$ and $id_1$. Then each arrow of ${{\mathcal C}_B}$ is a linear
combination of compositions of the form $ab^*$, with $a,b$ arrows in ${{\mathcal C}_B}^+$.
\end{lemm}

\begin{proof}
By definitions the arrows of ${{\mathcal C}_B}$ are linear
combinations of compositions of tensor products of maps
of the form $\eta$, $\mu^*$, $\eta^*$, $\mu$ and $id_1$. By an easy
induction argument, it suffices to prove is that each composition of the form $(id_?\otimes x\otimes
id_?)(id_?\otimes y\otimes id_?)$ with $x\in \{ \eta^*,\mu ,id\}$ and
$y\in \{ \eta,\mu^* ,id\}$ may be written as $\lambda (id_?\otimes z\otimes
id_?)(id_?\otimes t\otimes id_?)$ with $z\in  \{ \eta ,\mu^* ,id\}$ and
$t\in  \{ \eta^* ,\mu,id\}$ and $\lambda\in\c$ (i.e. that ``modulo
scalars, in each composition giving rise to an arrow of ${{\mathcal C}_B}$, the
$\eta^*$'s and $\mu$'s can be moved to the right''). All $3\times 3=9$
assertions to be verified are clear from lemma 2.1. 
\end{proof}

\begin{lemm}
For any $p\geq 0$ define $(\mu^*)^{(p)}\in Hom_{{\mathcal C}_B} (1,p+1)$ by
$(\mu^*)^{(0)}=id_1$ and by
$$(\mu^*)^{(p)}=(id_{p-1}\otimes\mu^* )(id_{p-2}\otimes\mu^* )\dots
(id\otimes\mu^* )\mu^*$$
if $p\geq 1$. Define also for any $p\geq 1$ an arrow $\eta^{(p)}\in
Hom_{{\mathcal C}_B} (0,p)$ by $\eta^{(p)}=(\mu^*)^{(p-1)}\eta$. Then for any $k\in\nn$ each
arrow in $Hom_{{\mathcal C}_B} (0,k)$ is a linear combination of compositions of
maps of the form $id_a\otimes \eta^{(p)}\otimes id_b$, with $a,b,p\in\nn$.
\end{lemm}

\begin{proof}
We know from lemma 3.1 that each arrow in $Hom_{{\mathcal C}_B} (0,k)$ is a linear
combination of arrows in ${{\mathcal C}_B}^+$. In particular each arrow in
$Hom_{{\mathcal C}_B} (0,k)$ is a linear combination of compositions of arrows of
the form $id_?\otimes \eta^{(?)}\otimes id_?$ and
$id_?\otimes\mu^*\otimes id_?$. It is enough to prove that each such
composition may be written without $\mu^*$'s. For, let us choose such
a composition, say $C$, having $z\geq 1$ $\mu^*$'s in its writing,
and assume that $C$ is not equal to a composition having $z-1$
$\mu^*$'s in its writing. Then $C$ is of the form $C^\prime
(id_a\otimes\mu^*\otimes id_b)D$ where $C^\prime$ has $z-1$
$\mu^*$'s in its writing, $a,b\in\nn$, and $D$ is a composition of
arrows of the form $id_?\otimes \eta^{(?)}\otimes id_?$. By using an easy
minimality argument on
the lenght of $D$ (``the term containing $\mu^*$ cannot be moved to the right'') we may assume that $D$ is of the form $(id_c\otimes
\eta^{(p)}\otimes id_d)D^\prime $ with $a+1\in\{ c+1,c+2,...,c+p\}$,
i.e. that
$$C=C^\prime (id_c\otimes ((id_{a-c}\otimes\mu^*\otimes
id_{b-d})(\mu^*)^{(p)}\eta )\otimes id_d)D^\prime$$
From the coassociativity property $(\mu^*\otimes
id)\mu^*=(id\otimes\mu^*)\mu^*$ of $\mu^*$ (cf. lemma 2.1 (iii)) we get
that the term in the middle is $(\mu^*)^{(p+1)}\eta =\eta^{(p+1)}$, contradiction.
\end{proof}

\begin{lemm}
Define a set $X_k\subset Hom_{{\mathcal C}_B} (0,k)$ for any $k\geq 1$ in the
following way. $X_1=\{\eta\}$ and for any $k\geq 2$
$$X_k=\{ (id_x\otimes\alpha\otimes id_y\otimes\beta\otimes
id_z\otimes\gamma\otimes\ldots )\eta^{(p)}\}$$
with $p$ ranging over $\{ 1,2,\ldots ,k\}$, $x,y,z,\ldots$ ranging
over sequences of strictly positive integers whose sum is $p$, and
$\alpha\in X_a$, $\beta\in X_b$, $\gamma\in X_c$ etc., with
$a,b,c,\ldots$ being positive integers whose sum is $k-p$.

Then $X_k$ is a system of generators of $Hom_{{\mathcal C}_B} (0,k)$ for any $k$.
\end{lemm}

\begin{proof}
We know from lemma 3.2 that for any $k$ the set $Y_k$ of compositions of
maps of the form $id_a\otimes\eta^{(p)}\otimes id_b$ with
$a,b,p\in\nn$ which happen to belong to $Hom_{{\mathcal C}_B} (0,k)$ is a system
of generators of $Hom_{{\mathcal C}_B} (0,k)$. We will prove by induction on $k$
that $X_k=Y_k$. For $k=1$ this is clear from lemma 3.2, so let $k\geq 2$. Pick an
arbitrary element of $Y_k$ 
$$C=\circ_{i=1}^{i=s}(id_{a_i}\otimes\eta^{(p_i)}\otimes id_{b_i})$$
We have $a_s=b_s=0$. It's easy to see that imposing the condition that
the biggest integer $t$ such that $a_t=0$ is minimal is the same as
assuming that $a_1,a_2,\ldots ,a_{s-1}$ are strictly positive
numbers. With $x:=inf\{ a_1,a_2,\ldots ,a_{s-1}\}$ and $p:=p_s$ we get that
$$C=(id_x\otimes f)\eta^{(p)}$$
with $f$ being a composition of maps of
the form $id_a\otimes\eta^{(q)}\otimes id_b$ which belongs to
$Hom_{{\mathcal C}_B} (p,k-x)$, i.e. with $f$ being of the form $\alpha\otimes
id_y\otimes\beta\otimes id_z\otimes\gamma\otimes\ldots$ with
$x,y,z,\ldots$, $\alpha ,\beta ,\gamma ,\ldots$ and $a,b,c,\ldots$
being as in the definition of $X_k$'s.
\end{proof}

\begin{lemm}
$dim(Hom_{{\mathcal C}_B} (0,k))\leq C_k :=\frac{(2k)!}{k!(k+1)!}$ (the $k$-th
Catalan number) $\forall\, k$.
\end{lemm}

\begin{proof}
With $D_k=dim(Hom_{{\mathcal C}_B} (0,k))$ lemma 3.3 shows that
$$D_k\leq
D_{k-1}+\sum_{a+b=k-2}D_aD_b+\sum_{a+b+c=k-2}D_aD_bD_c+\cdots$$
(here each sum corresponds to a value of $p$ in lemma 3.3, and
$a,b,c,\ldots$ are the ones in lemma 3.3). Thus $D_k\leq E_k$, where $E_k$ are the numbers defined by $E_0=E_1=1$
and by the above formula with $E_m$ at the place of $D_m$ and with $=$
at the place of $\leq$. By rearranging terms in this equality we get
(by an easy induction on $s$) that
$$E_s=\sum_{x+y=s-1}E_xE_y$$
for any $s$. It follows that the $E_s$'s are the Catalan numbers (well-known, just consider the square of the series $\sum E_sz^s$...) and we are done.
\end{proof}

\begin{prop}
Let $k\in\nn$.

(i) The inclusion $A_{n,2k}\subset End_{{\mathcal C}_B} (k)$ given by proposition 2.2 is
an equality.

(ii) The set $X_k$ in lemma 3.3 is a basis of $Hom_{{\mathcal C}_B} (0,k)$.

(iii) $dim(Hom_{{\mathcal C}_B} (0,k))=C_k$.
\end{prop}

\begin{proof}
For any $l\in\nn$ we have
$$dim(End_{{\mathcal C}_B} (l))=dim(Hom_{{\mathcal C}_B} (0,2l))\leq C_{2l}\leq
dim(End_{{\mathcal C}_B} (l))$$

Indeed, by Frobenius reciprocity we have $End_{{\mathcal C}_B} (l)\simeq Hom_{{\mathcal C}_B}
(0,2l)$ for any $l$, and this gives the equality on the left. Lemma 3.4 and proposition 2.2 give the inequalities.

Thus both inequalities have to be equalities, and this proves (i). This
proves also (ii,iii) for even values of $k$. By definition of $X_k$ we
have $X_k\otimes\eta\subset X_{k+1}$, and with $k=2l+1$ we get that
the set $X_{2l+1}\otimes\eta$ consists of linearly
independent maps. Thus  $X_{2l+1}$ consists of linearly
independent maps, so we get (ii,iii) for odd values of $k$.
\end{proof}

\section{Corepresentations of $A^{aut}(B)$. Applications}

Let $B$ be a finite dimensional $\c^*$-algebra of dimension $\geq 4$.

\begin{theo}
The set of classes of finite dimensional irreducible smooth
corepresentations of $A^{aut}(B)$ can be labeled by the positive
integers, $Irr(A^{aut}(B))=\{ p_n\mid n\in\nn\}$, such that the fusion
rules are
$$p_k\otimes p_s=p_{\mid k-s\mid}+p_{\mid k-s\mid
  +1}+\dots +p_{k+s-1}+p_{k+s}$$
\end{theo}

\begin{proof}
We will use proposition 3.1 (iii) and a method from \cite{cras}.

We first recall a few well-known facts on $\ss\oo (3)$. Let
$\chi_0,\chi_1,\chi_2,\chi_3,\ldots$ be the characters of the
irreducible representations of $\ss\oo (3)$, listed in the increasing order of dimensions. They satisfy the formulas
in theorem 4.1, with $\chi_m$ at the place of $p_m$ and with $\cdot$
at the place of $\otimes$. We denote by $C(\ss\oo (3))_c$ the
$*$-subalgebra of $C(\ss\oo (3))$ generated by the $\chi_i$'s. Let
$\int$ be the integration over $\ss\oo (3)$. Then the $\chi_i$'s form an
orthonormal basis of $(C(\ss\oo (3))_c,\int )$. There exists a
canonical isomorphism
$$\c [X]\simeq C(\ss\oo (3))_c,\,\,\, X\mapsto\chi_1$$
For any $k$ the multiplicity of $1$ into $(1+\chi_1)^k$ is $C_k$. See
e.g. \cite{cras} with $F=\begin{pmatrix}0&1\cr
  -1&0\end{pmatrix}$ for everything (we didn't succeed in finding a classical reference).

Let $A^{aut}(B)_c$ be the algebra of characters of
corepresentations of $A^{aut}(B)$, let $h:A^{aut}(B)\to\c$ be the Haar
functional, and denote by $r\to\chi (r)$ the
character of corepresentations (see \cite{w1}). By using the above
isomorphism, we may define a morphism of algebras in the following way.
$$\phi : C(\ss\oo (3))_c\to A^{aut}(B)_c,\,\,\, \chi_1\mapsto \chi (u)-1$$
We have on one hand that the multiplicity of $1$ into $(1+\chi_1)^k$
is $C_k$, and on the other hand that the multiplicity of $1$ into
$u^{\otimes k}$ is also $C_k$ (cf. proposition 3.1 (iii) and the
definition of ${{\mathcal C}_B}$). This could be interpreted as saying that
$$h\phi ((\chi_1+1)^k)=\int (\chi_1+1)^k$$
for any $k$. It's easy to get from this (by induction on $s$) that $h\phi
(\chi_1^s)=\int (\chi_1^s)$ for any $s$, and as $\chi_1$ generates $C(\ss\oo (3))_c$ as
an algebra we get that
$$h\phi =\int$$
Thus $\{\phi (\chi_k)\mid k\geq 0\}$ is an orthonormal
basis of $A^{aut}(B)_c$, and for finishing the proof it's enough to
construct corepresentations $p_k$ of $A^{aut}(B)$ such that $\chi
(p_k)=\phi (\chi_k)$ for any $k$. Indeed, the fact that
$Irr(A^{aut}(B))=\{ p_n\mid n\in\nn\}$ will be clear from this and from
Peter-Weyl type theory from \cite{w1}; and the assertion on fusion
rules will be also clear from this
and from the fusion rules for irreducible representations of $\ss\oo (3)$.

We do it by induction on $k$. We may define $p_0=1$. Also as $\eta\in
Hom(1,u)$ we know that $u$ contains a copy of $1$, so we may set
$p_1=u-1$. So let $k\geq 2$ and assume that we have constructed $p_0,p_1,\ldots ,p_{k-1}$ with
$$\chi (p_i)=\phi (\chi_i),\,\,\, i=0,1,\ldots ,k-1$$
We know from fusion rules for representations of $\ss\oo (3)$ that
$$\chi_{k-1}\chi_1=\chi_{k-2}+\chi_{k-1}+\chi_k$$
so by applying $\phi$ and by using $h\phi =\int$ we get that both
$p_{k-2}$ and $p_{k-1}$ are subcorepresentations of $p_{k-1}\otimes
p_1$. Thus there exists a corepresentation $p_k$ such that
$$p_{k-1}\otimes p_1=p_{k-2}+p_{k-1}+p_k$$
From the above two formulas we get $\chi (p_k)=\phi (\chi_k )$ and we are
done.
\end{proof}

We recall from Wang \cite{wn1} that for $m\geq 2$ the $\c^*$-algebra
$A_o(m)$ is defined with self-adjoint generators $\{
v_{ij}\}_{i,j=1,...,m}$ and the relations making the matrix
$v=(v_{ij})$ unitary. The pair $(A_o(m),v)$ is a finitely generated
full Woronowicz-Kac algebra.

\begin{coro}
For any $m\geq 2$ the enveloping $\c^*$-algebra of the $*$-subalgebra
of $A_o(m)$ generated by the coefficients of the square of its
fundamental corepresentation is canonically isomorphic (as a full Woronowicz-Kac algebra) to $A^{aut}(M_m(\c))$.
\end{coro} 

\begin{proof}
Let $\x =*-alg(\{ v_{ij}v_{kl}\}_{i,j,k,l=1,...,m})$ be the
$*$-algebra in the statement and let $X$ be its enveloping
$\c^*$-algebra. Then the pair $(X,v^{\otimes 2})$ is a finitely
generated full Woronowicz-Kac algebra. We have the following two facts.

I. We recall that the fusion semiring $R^+(A)$ of a Woronowicz-Kac algebra $A$
is the set of equivalence classes of corepresentations of $A$ together
with the operations sum and tensor product. By \cite{cras} we have an
isomorphism
$$R^+(A_o(m))\simeq R^+(C(\ss\uu
(2)))$$
sending $v$ onto the fundamental corepresentation of $C(\ss\uu
(2))$. It follows that we have an isomorphism
$R^+(X)\simeq R^+(C(\ss\oo (3)))$ sending $v^{\otimes 2}$ onto the sum
of the fundamental corepresentation of $C(\ss\oo (3))$ with the
trivial corepresentation. On the other hand, theorem 4.1 gives an isomorphism
$$R^+(A^{aut}(M_m(\c)))\simeq R^+(C(\ss\oo (3)))$$
which sends $u$
onto the sum of the fundamental corepresentation of $C(\ss\oo (3))$ with the
trivial corepresentation. By combining these two facts we get an isomorphism
$$g:R^+(A^{aut}(M_m(\c)))\simeq R^+(X)$$
which sends $u\mapsto v^{\otimes 2}$.

II. Consider the fundamental coaction $\phi$ of $A_o(m)$ on
$M_m(\c )$, i.e. the map $\phi :M_m(\c )\to M_m(\c )\otimes A_o(m)$ given by
$x\mapsto ad(v)(x\otimes 1)$. The image of $\phi$ is contained in
$M_m(\c )\otimes\x$, so we get by restriction a coaction $\psi :M_m(\c
)\to M_m(\c )\otimes X$. It follows that there exists a morphism of $\c^*$-algebras
$$f:A^{aut}(M_m(\c))\to X$$
such that $(id\otimes f)\beta =\psi$, where $\beta$ is the fundamental
coaction of $A^{aut}(M_m(\c))$ on $M_m(\c )$. Thus at the level of
corepresentations we get a morphism of semirings
$$f_*:R^+(A^{aut}(M_m(\c)))\to R^+(X)$$
which sends $u\mapsto v^{\otimes 2}$.

By using I and II we get that $g^{-1}f_*$ is an endomorphism of
$R^+(A^{aut}(M_m(\c)))$ sending $u\mapsto u$. By using the structure
of $R^+(A^{aut}(M_m(\c)))$ given by the fusion rules in theorem 4.1 we
get from this that $g^{-1}f_*=id$, so in particular $f_*$ is an
isomorphism. It follows that $f$ maps any basis of $\a^{aut}(M_m(\c))$
consisting of coefficients of irreducible corepresentations
(cf. Woronowicz' Peter-Weyl theory, see \cite{w1}) onto a basis of $\x$
consisting of coefficients of irreducible corepresentations. Thus the
restriction $\a^{aut}(M_m(\c))\to\x$ of $f$ is an isomorphism, and by
passing at the enveloping $\c^*$-algebras level we get that $f$ itself is an
isomorphism.
\end{proof}

\begin{coro}
$A^{aut}(B)$ is amenable as a Woronowicz algebra iff $dim(B)=4$.
\end{coro}

\begin{proof}
Theorem 4.1 gives an isomorphism
$$R^+(A^{aut}(B))\simeq R^+(C(\ss\oo (3)))$$
which sends the fundamental $dim(B)$-dimensional corepresentation of
$A^{aut}(B)$ onto a $4$-dimensional corepresentation of $C(\ss\oo
(3))$ (see the proof of corollary 4.1). We may conclude by applying the following consequence of the quantum Kesten result
 (proposition 6.1 in \cite{subf}): is $A,B$ are two finitely generated
 Woronowicz algebras, if $A$ is amenable and if $\varphi : R^+(A)\simeq
 R^+(B)$ is an isomorphism, then $B$ is amenable if and only if $\varphi$ is
 dimension-preserving.
\end{proof}

\begin{coro}
Let $P$ be a $II_1$ factor, $\alpha$ be a minimal coaction of
$A^{aut}(B)$ on $P$ and $\beta$ be the canonical coaction of
$A^{aut}(B)$ on $B$. Then the fixed point subfactor (see \cite{spin})
$$P^\alpha\subset (B^o\otimes P)^{\beta\odot\alpha}$$
has index $dim(B)$ and principal graph $A_\infty$.
\end{coro}

\begin{proof}
This is clear from theorem 4.3 (ii) in \cite{spin} which asserts that
the relative commutants of even order of any fixed point subfactor of the form
$P^\alpha\subset (B^o\otimes P)^{\beta\odot\alpha}$ are the algebras
of intertwiners of the tensor powers of the corepresentation
corresponding to the coaction $\beta$. In our present situation, the
corepresentation corresponding to $\beta$ is the fundamental
corepresentation $u$ of $A^{aut}(B)$, the algebras
of intertwiners of its tensor powers are Temperley-Lieb algebras
(cf. proposition 3.1 (i); in fact this is also clear from theorem
4.1), so our subfactor has graph $A_\infty$.
\end{proof}

\end{document}